\RequirePackage{fix-cm}
\documentclass[twocolumn]{svjour3}          
\smartqed  
\usepackage{graphicx}
\usepackage{float}
 \usepackage{mathptmx}      
\usepackage{latexsym}
 \journalname{Nonlinear Dynamics}
\begin{document}

\title{On the use of the theory of dynamical systems for transient problems 
}
\subtitle{A preliminary work on a simple model}


\author{Ugo Galvanetto          \and
       Luca Magri 
}


\institute{Ugo Galvanetto \at
              Dipartimento di Ingegneria Industriale, Universit\`a di Padova \\
              Via Marzolo 9, 35131 Padova, Italy \\
              Tel.: +39 049 8275594\\
              Fax: +39 049 8275604\\
              \email{ugo.galvanetto@unipd.it}           
           \and
           Luca Magri \at
              Department of Engineering, University of Cambridge\\
              Trumpington Street, CB2 1PZ, Cambridge, U.K.\\
              Tel.: +44 (0)1223 746971\\
              Fax:  	+44 (0)1223 765311\\
              \email{lm547@cam.ac.uk}
}

\date{Received: date / Accepted: date}

\maketitle

\begin{abstract}
This paper is a preliminary work to address the problem of dynamical systems with parameters varying in time. An idea to predict their behaviour is proposed. 
These systems are called \emph{transient systems}, and are distinguished from \emph{steady systems}, in which parameters are constant. 
In particular, in steady systems the excitation is either constant (e.g. nought) or periodic with amplitude, frequency and phase angle which do not vary in time. 
We apply our method to systems which are subjected  to a transient excitation, which is neither constant nor periodic. 
The effect of switching-off and full-transient forces is investigated. 
The former can be representative of switching-off procedures in machines; the latter can represent earthquake vibrations, wind gusts, etc. acting on a mechanical system. 
This class of transient systems can be seen as the evolution of an ordinary steady system into another ordinary steady system, for both of which the classical theory of dynamical systems holds. 
The evolution from a steady system to the other is driven by a transient force, which is regarded as a map between the two steady systems. 
\end{abstract}
\section{Introduction}
\label{intro}
Dynamical systems are extremely pervasive in science and technology and are used to model natural and engineering phenomena. Very often mathematical models are autonomous, i.e. they do not entail an explicit dependence on time. Such a definition can be extended to include as well systems with periodic excitation, since time can be added as an extra-variable with constant time-derivative  \cite{STROGATZ,GUCK}. 
In dissipative dynamical systems, which are of great use in engineering, it is possible to introduce a qualitative distinction between \emph{steady states} and \emph{transient states}.
A steady state is characterised by recurrent behaviour, i.e. a particular point in the phase space is a steady state if the system, after sufficient time, returns arbitrarily close to the point \cite{THOMPSON}. This definition includes fixed points, limit cycles, quasi-periodic and chaotic steady states. A point in the phase space is a transient state if it is not a steady state. 
Both steady and transient states are typical of dissipative dynamical systems with time-invariant parameters and constant or periodic excitations, i.e. autonomous or periodically forced systems which therefore can be called \emph{steady systems}. 

Transient systems, as opposed to steady systems, are the topic of the present work. They are characterised by system and/or excitation parameters which change in time. In particular, we limit our attention to systems which are transient because they are subjected to non-periodic excitations. 
The theory of dynamical systems is well suited to deal with steady problems. Fixed points, limit cycles, basins of attraction ... are all features of steady problems, even if they can be used to study transient trajectories. 
However, in the real world engineers are usually faced by a large variety of problems, some of which can be described as steady, whereas many others as transient.

This paper proposes an  idea about how to make use of concepts of the theory of dynamical systems to address transient problems.

The paper is organised as follows: section~\ref{types} presents the main concepts of the paper, section~\ref{simple} proposes a prototypical model which will be subjected to transient forces in the examples of section~\ref{numerical}. A brief discussion and conclusions will end the paper in section~\ref{concl}.
%
\section{Examples of transient forces} \label{types}
%
We consider engineering systems which can be modelled with nonlinear ordinary differential equations (ODEs), such as:
\begin{equation}\label{eq1}
\frac{\mathrm{d}x}{\mathrm{d} t}=f(x,\gamma),                           
\end{equation}
where $x$ is the state vector, $\gamma$ is the set of system parameters, and $t$ is  time.
If the system is forced then an additional term appears in eq.~\ref{eq1}:
\begin{equation}\label{eq2}
\frac{\mathrm{d}x}{\mathrm{d} t}=f(x,\gamma)  + a\cdot p(t),      
\end{equation}
where $a$ is the excitation amplitude and $p(t)$ is a normalised periodic function of time such that $p(t)=p(t+T)$, where $T$ is the period. 
If both the system and excitation parameters 
are constant, the system is called a \emph{steady system}. 
On the other hand, if the system parameters and/or the excitation parameters change in time, the system described by either eq.~(\ref{eq1}) or ({\ref{eq2}) is a \emph{transient system}. 
%

In this paper, we will assume that the system can operate in two steady conditions: one with no excitation, represented by eq.~(\ref{eq1}), and the second under the action of a periodic excitation, represented by eq.~(\ref{eq2}). 
We will assume as well that in both conditions there exists  only one \emph{acceptable functioning mode}. 
For example, an acceptable functioning mode can be 
an attracting fixed point, when the excitation is not present, 
or a stable limit cycle, when the periodic excitation is present. 
Moreover, in the unforced system a \emph{safe zone} exists around the stable fixed point, within its basin of attraction. 
The safe zone represents all configurations that the system can safely assume when in action and affected by perturbations or imperfections that keep it slightly away from the stable fixed point.
Likewise, for a periodically forced system a safe zone exists around the stable limit cycle, within its basin of attraction \cite{SOLIMAN,LENCI}. 
%
We assume that safe zones are limited sets so that they are contained by boundaries of finite length (unless the basin boundaries are fractal, case which is not considered in this paper).

The main idea of this work is that a transient problem can be regarded as a map between two steady systems: 
one is the initial system  and the other the final system. 
We apply this idea to a case of transient problem in which the system parameters are constant and only the excitation amplitude, $a(t)$, can vary. 
The periodic function $p(t)$ is chosen to be a harmonic function with constant angular frequency, $p(t)=\cos(\omega t)$.
Hence, the transition between the two steady systems is driven by a transient force with a known time-varying amplitude (stochastic forces are not studied because beyond the scope of this  work).
We consider two different types of transient forces (fig.~\ref{fig:1}).

Firstly, a \emph{switching-off} transient force brings a forced system to an unforced condition in a finite time. 
During this time, the amplitude of the transient force varies from the value $a_0\not=0$, at the beginning of the transient, 
to zero, at the end of the transient as shown in fig.~\ref{fig:1}a. The initial system is represented by eq.~\ref{eq2} and the final system by eq.~\ref{eq1}. (Similarly, a switching-on transient force can be defined by an increasing force amplitude.)

Secondly, a \emph{full}-transient force acts between two unforced systems for a finite time. 
During this time, the amplitude of the transient force varies continuously assuming the value 0 both at the beginning and the end of the transient, 
as depicted in the sketch of fig.~\ref{fig:1}b. Both initial and final systems are given by eq.~\ref{eq1}.
Similarly, the full-transient force can act also between two forced systems. 

The above two types of transient forces are representative of common engineering transient problems. 
\begin{figure*}[]
\centering
  \includegraphics[width=0.75\textwidth]{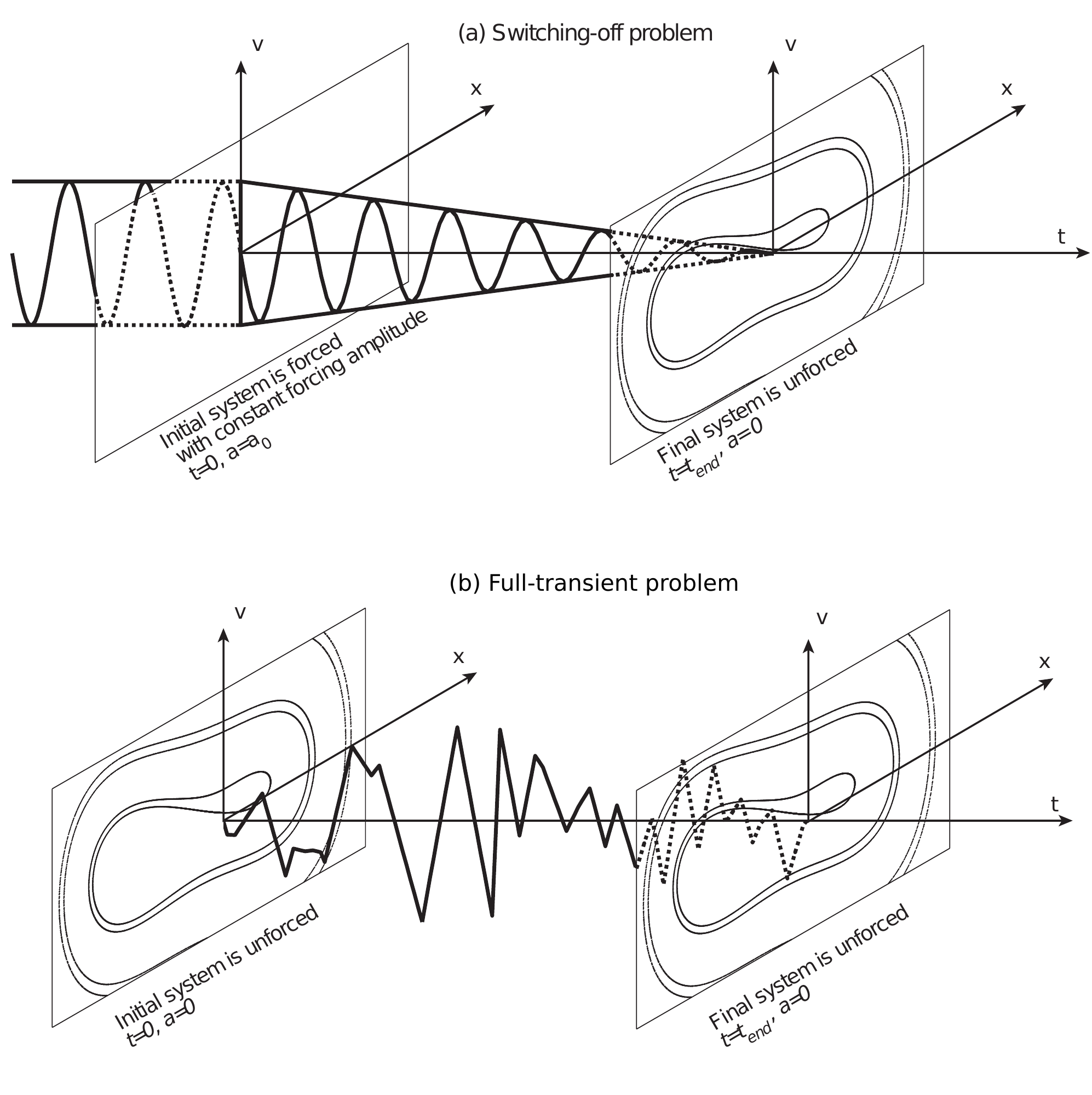}
\caption{Sketch of transient systems with (a) switching-off transient force and (b) full-transient force. The stable manifolds of a saddle point are drawn in the phase plane of the unforced system.}
\label{fig:1}       
\end{figure*}
\begin{figure*}[]
\centering
  \includegraphics[width=0.85\textwidth]{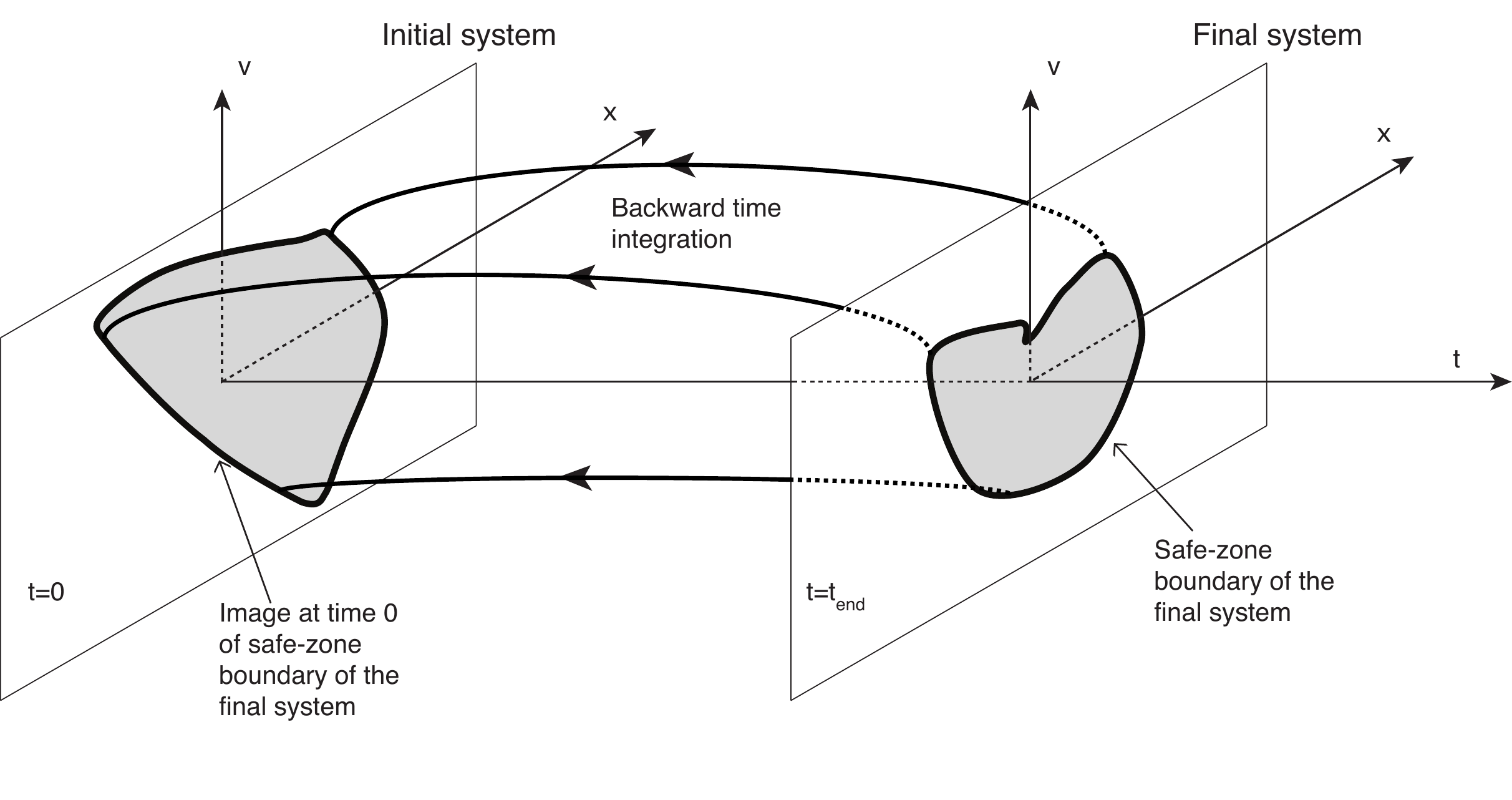}
\caption{The boundary of the safe-zone at $t_{end}$ is integrated backward in time to find its image at $t=0$.}
\label{fig:2}       
\end{figure*}
\begin{figure*}[]
\centering
  \includegraphics[width=0.5\textwidth]{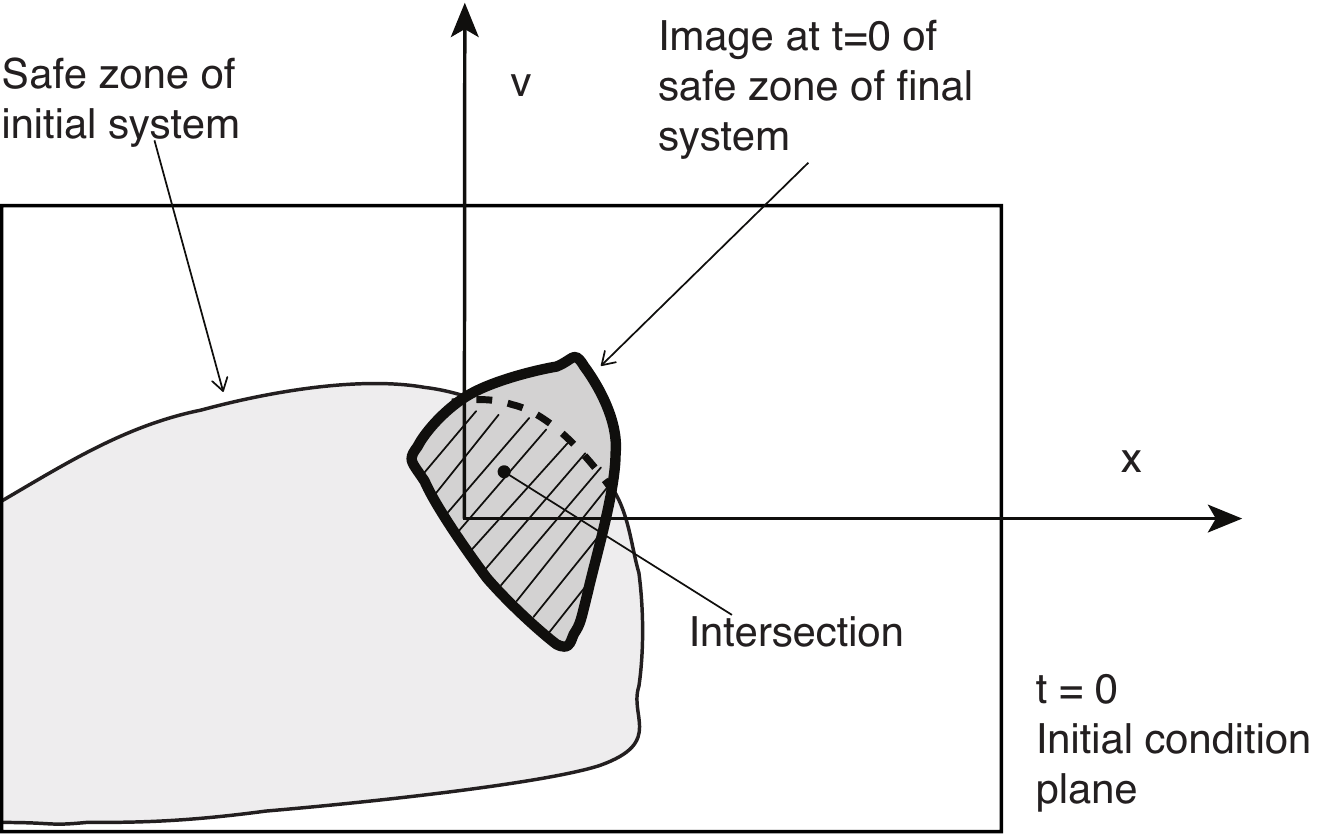}
\caption{The image at $t=0$ of the  safe-zone at $t_{end}$ may overlap the safe zone of the initial system only in part. If the intersection is empty, then it is not possible to bring the system from the initial safe zone to the final one under the action of the considered transient force. }
\label{fig:3}       
\end{figure*}
The transient problem is characterised by a finite duration in time and it is bound by two steady systems. 
For both types of transient forces above-mentioned, the theory of dynamical systems can be applied to the two steady systems, which exist before and after the transient, so that we can have a complete description of the dynamical features of the steady systems: attractors, basins of attraction, basin boundaries, etc and therefore of the relevant safe zones. 
In many engineering applications, however, it is essential to determine if a transient between two steady conditions brings the system to an acceptable functioning mode or not. 
Thus, the engineering point of view of a transient problem can be described in the following way: 
can we determine the set of initial conditions (at $t=0$) which, under the action of the transient force, give origin to a trajectory ending up in the safe zone of the final system (for $t>t_{end}$, where $t_{end}$ is the end of the transient duration)? 
The problem is illustrated in the sketch of fig.~\ref{fig:2}. 
%
In general, we can locate the steady states of the final system \cite{FOALE,SECONDO}, in particular its saddles (points or cycles)  with the relevant stable manifolds. 
One of the stable steady states is acceptable and its safe region will therefore be inside its basin of attraction. 
To answer the afore-mentioned question, we have to integrate backward in time all points belonging to the boundary of the safe zone from $t_{end}$ to $t=0$ 
under the action of the time-varying-amplitude (transient) force.  
The image at $t=0$ of the boundary of the safe zone of the final system (called \emph{transformed boundary} in the remainder of the paper) 
 will provide the boundary, in the initial-condition plane ($t=0$), of the set of safe points, called \emph{safe transient initial conditions}. 
These are the points of the initial system that, if the transient force is activated, give origin to trajectories ending up in the safe region of the final system. It is particularly important to find if there exists an intersection between the safe zone of the initial system, and the image at $t=0$ of the safe zone of the final system (obviously contained in the transformed boundary), as shown in fig.~\ref{fig:3}. If such an intersection is empty then it is not possible to transfer the system from the initial safe zone to the final one under the action of the transient force.
%
\section{A simple model as a prototype}
\label{simple}
In the present section we introduce a simple model as a prototype which satisfies all assumptions made in \S~\ref{types}. The main idea presented in \S~\ref{types} will be applied to the prototype in \S~\ref{numerical}. The prototype is a nonlinear oscillator, chosen for its simplicity, 
and it is characterised by (at least) two coexisting attractors, one of which is representing an acceptable solution whereas the other is a solution which is not acceptable from an engineering point of view. 

The system is governed by the following equation:
\begin{equation} \label{eq:123}
m\frac{\mathrm{d}^2 x}{\mathrm{d}t^2} + c\frac{\mathrm{d} x}{\mathrm{d}t} - k_1x+k_2x^2 + k_3x^3=a\cdot p(t),
\end{equation}
where $a\cdot p(t)$ is a function of time and represents the forcing term which may contain the transient force. 
The system parameters have the following values: $m=1$, $c=0.25$, $k_1= k_2=k_3=1$.
The unforced nonlinear oscillator has three fixed points 
whose coordinates are $x_1=0$, $x_2=(-1+\sqrt{5})/2\approx 0.618$ and $x_3=(-1-\sqrt{5})/2\approx -1.618$.
The fixed point at the origin is a saddle point and $x_{2,3}$ are two point attractors. 
The stable manifolds of the saddle point divide the phase plane of the unforced system 
in the two basins of attraction of the fixed points $x_{2,3}$, as shown in fig.~\ref{fig:bas}. 
\begin{figure}[]
\centering
  \includegraphics[width=0.5\textwidth]{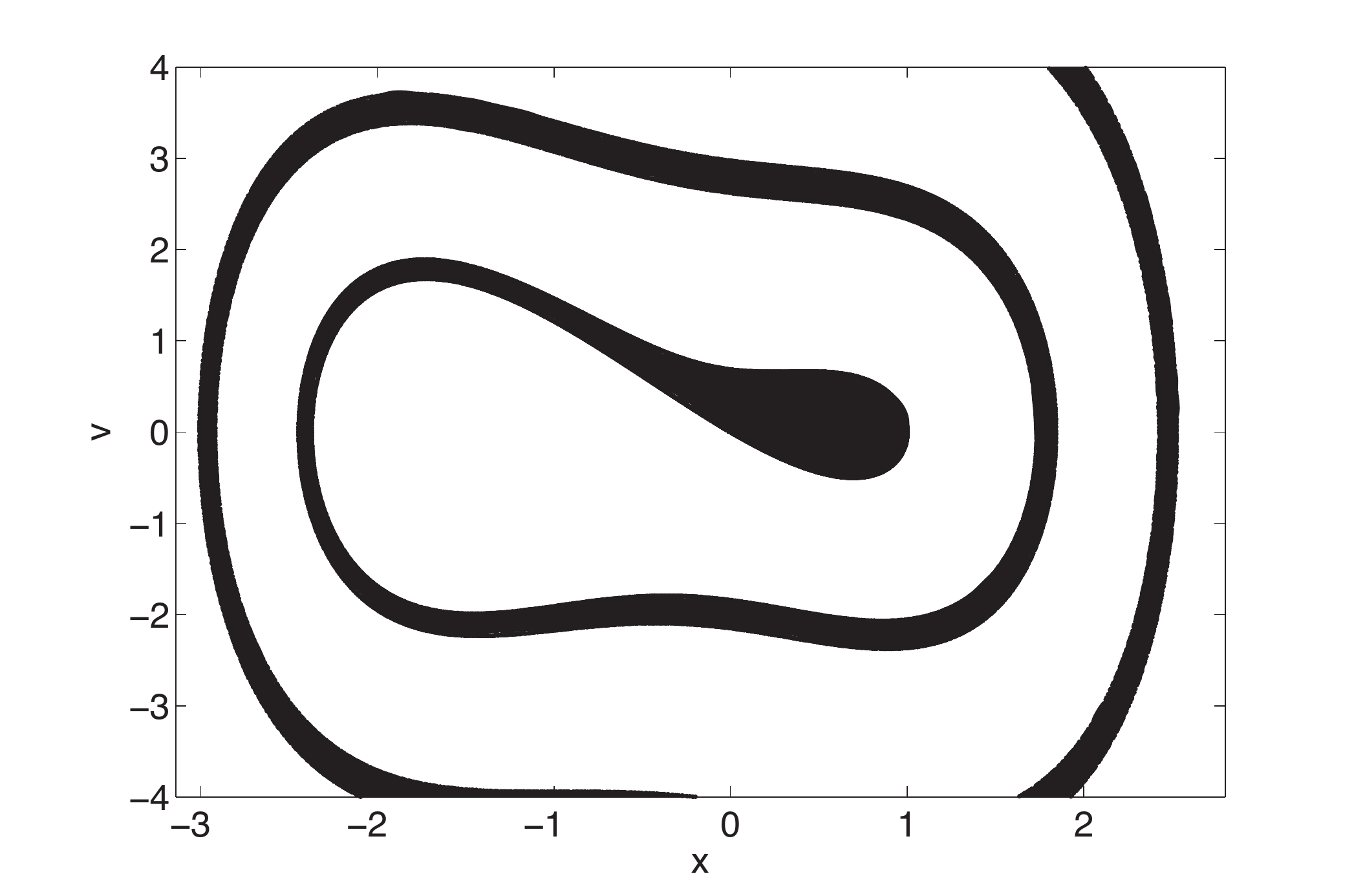}
\caption{Basins of attraction of the  attracting fixed points: the black region is the basin of $x_2$ and the white region is the basin of $x_3$.}
\label{fig:bas}       
\end{figure}
If a small-amplitude harmonic external force is added to the system the fixed points are transformed in limit cycles of the same nature: 
the saddle point becomes a saddle limit cycle and the two point attractors become two attracting limit cycles. 
The stable manifolds of the saddle cycle divide the phase space of the forced system in the two basins of attraction of the attracting limit cycles. A Poincar\'e section of the phase space would provide a figure similar to that of the basins of attraction of the unforced system, so that fig.~\ref{fig:bas} is a good approximation of the basins of attractions of the forced system.

It is assumed that the fixed point $x_2=(-1+\sqrt{5})/2$ and the corresponding stable limit cycle represent the acceptable functioning modes of the system, 
whereas the fixed point $x_3=-(1+\sqrt{5})/2$ and the corresponding stable limit cycle represent the unacceptable functioning modes.
The system is \emph{off} when the external force $a\cdot p(t)$ is nought.
The system is \emph{on} if $a\cdot p(t)$ is an harmonic force with constant amplitude,
$a\cdot p(t) = a_0 \cos(\omega t)$. 

The transition between the states off and on of the system takes place with a continuous variation in time of the amplitude $a(t)$, 
which is a function of time only during the transient and is constant before and after it. 

The switching-off problem (fig.~\ref{fig:1}a) is defined as follows:
\begin{enumerate}
\item the initial system is forced with constant amplitude, $a(t)=a_0$, when $t<0$;
\item $a(t)=a_0\left(1-\frac{t}{t_{end}}\right)$ when $0 \leq t \leq t_{end}$;
\item the final system is unforced, $a(t)=0$ when $t>t_{end}$;
\end{enumerate} 
where 
 $0 \leq t \leq t_{end}$ represents the duration of the transient. 
We chose the force amplitude to decrease linearly, but any other continuous variation could be assumed. 
The full-transient problem (fig.~\ref{fig:1}b) is defined as follows:
\begin{enumerate}
\item the initial system is unforced, $a(t)=0$, when $t<0$;
\item the force varies continuously when $0 \leq t \leq t_{end}$;
\item the final system is unforced when $t>t_{end}$.
\end{enumerate}
%
\section{Results and discussion} \label{numerical}
The prototype model presented in \S~\ref{simple} is subjected to the transient forces which represent two examples of  transient systems described in \S~\ref{types}. These types of forces have important roles in many engineering applications. However, the method we propose can be applied conceptually in the same manner to other types of transient systems (see \S~\ref{types}), once one knows: 
i) the duration of the transient;
ii) the mathematical law that describes the transient state;
iii) the initial and final steady systems.

The 
force parameters are $\omega=3$ $rad/s$ and $a_0= 0.2$.  In all examples of the present section the safe region of the initial system is chosen to coincide with the basin of attraction of the relevant acceptable attractor.

The results presented in this section are obtained by numerical time integration, which is performed via a $4^{th}-5^{th}$ order Runge-Kutta-Fehlberg algorithm with adaptive step  \cite{ODE45}. 
The relative and absolute integration errors are taken to be $\sim O(10^{-8})$.

The proposed procedure is based on a reverse-time integration which is well defined only if the mechanical systems under investigation are well defined in negative time; therefore it cannot be immediately applied to non-smooth systems such as those affected by dry friction \cite{Galvanetto1,Galvanetto2}.
\subsection{Switching-off transient}\label{sot}
At the beginning of the transient the forcing amplitude is $a_0$ and it decreases linearly in time reaching the value zero after $t_{end}(=2.0)$  seconds. 
We want to know the set of initial conditions in the phase plane at $t=0$ from which it is safe to switch off the excitation. These are the safe transient  initial conditions, as defined in \S~\ref{types}. 
An accurate approximation of the stable manifold of the unforced system (final system) is found with standard computational techniques \cite{LENCI,FOALE} and it is shown in fig.~\ref{fig:bas}. The safe zone around the fixed point in the unforced phase plane is arbitrarily chosen coincident with the basin itself capped with the line $v=1$, as shown in fig.~\ref{fig:OFFbas}a. The required set of safe transient initial conditions is the set of points in the phase plane at $t=0$ which, under the action of the decreasing excitation, are origin of trajectories ending up within the safe zone of the unforced system. 
The boundary of the safe zone of the unforced system is shown as a black line in figure \ref{fig:OFFbas}a. 
The integration backward in time of the boundary of the safe region is performed choosing as initial conditions a set of closely spaced points along the black curve in fig.~\ref{fig:OFFbas}a and integrating the following system:
\begin{eqnarray}
\label{eq:dotv}
&\dot{x}=   - v, \\
&\dot{v}=   - \left(a(t)\cos(\omega t) - c v + x - x^2 - x^3)\right),\\ 
&a(t)=a_0\left(\frac{t}{t_{end}}\right),\\
&0 \leq t \leq t_{end},\;t_{end}=2\pi/\omega.
\end{eqnarray}
The integration of the above equations provides the transformed boundary shown as a black curve in fig.~\ref{fig:OFFbas}b.

\begin{figure}
\centering
  \includegraphics[scale=0.45]{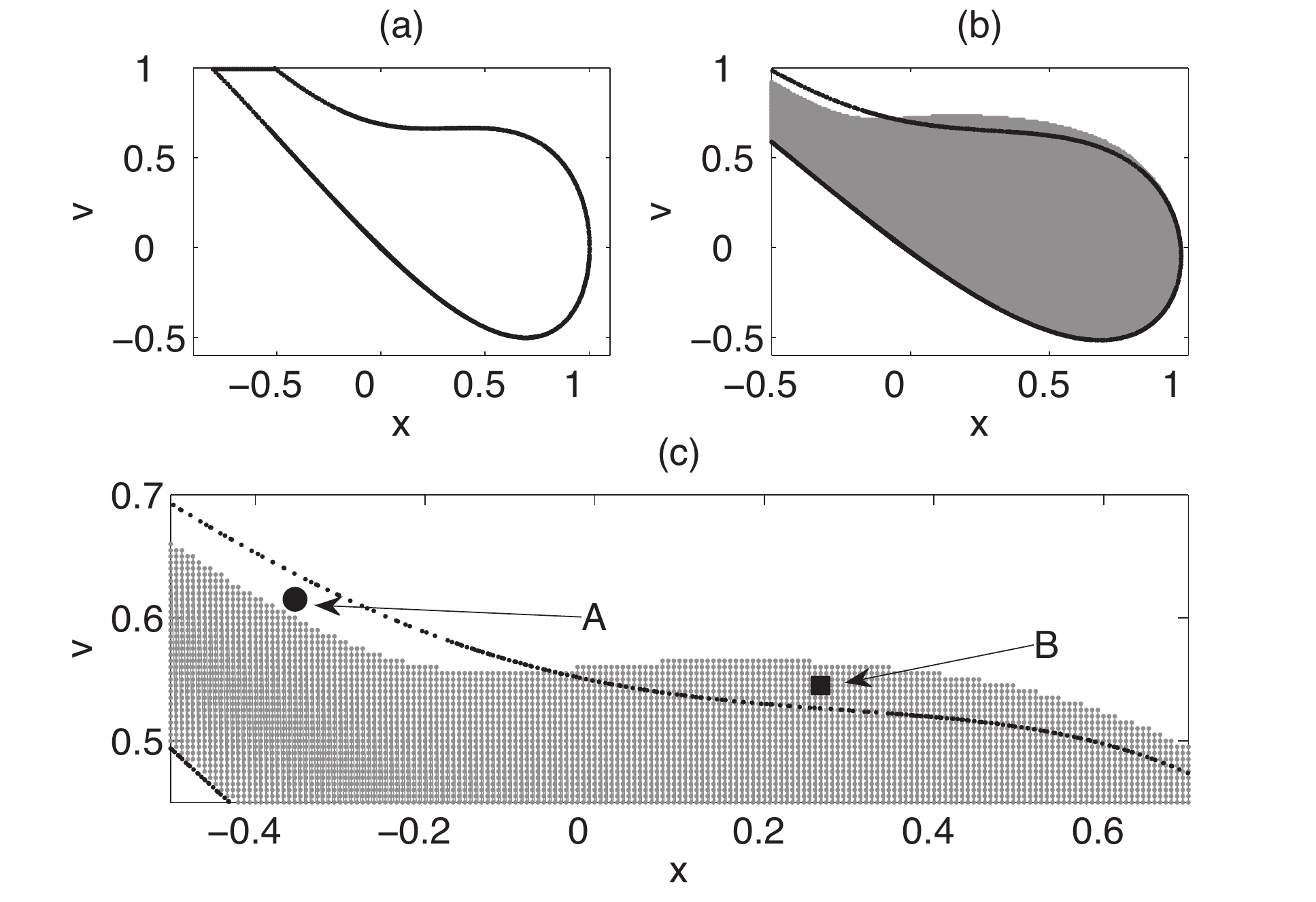}
\caption{Switching-off case: (a) safe zone in the final system; (b) superposition of the transformed boundary of the safe zone of the final system (dotted line) with the basin of attraction of the initial system (grey zone); (c) magnified view of a portion of frame (b).}
\label{fig:OFFbas}       
\end{figure}
\begin{figure}
  \includegraphics[scale=0.425]{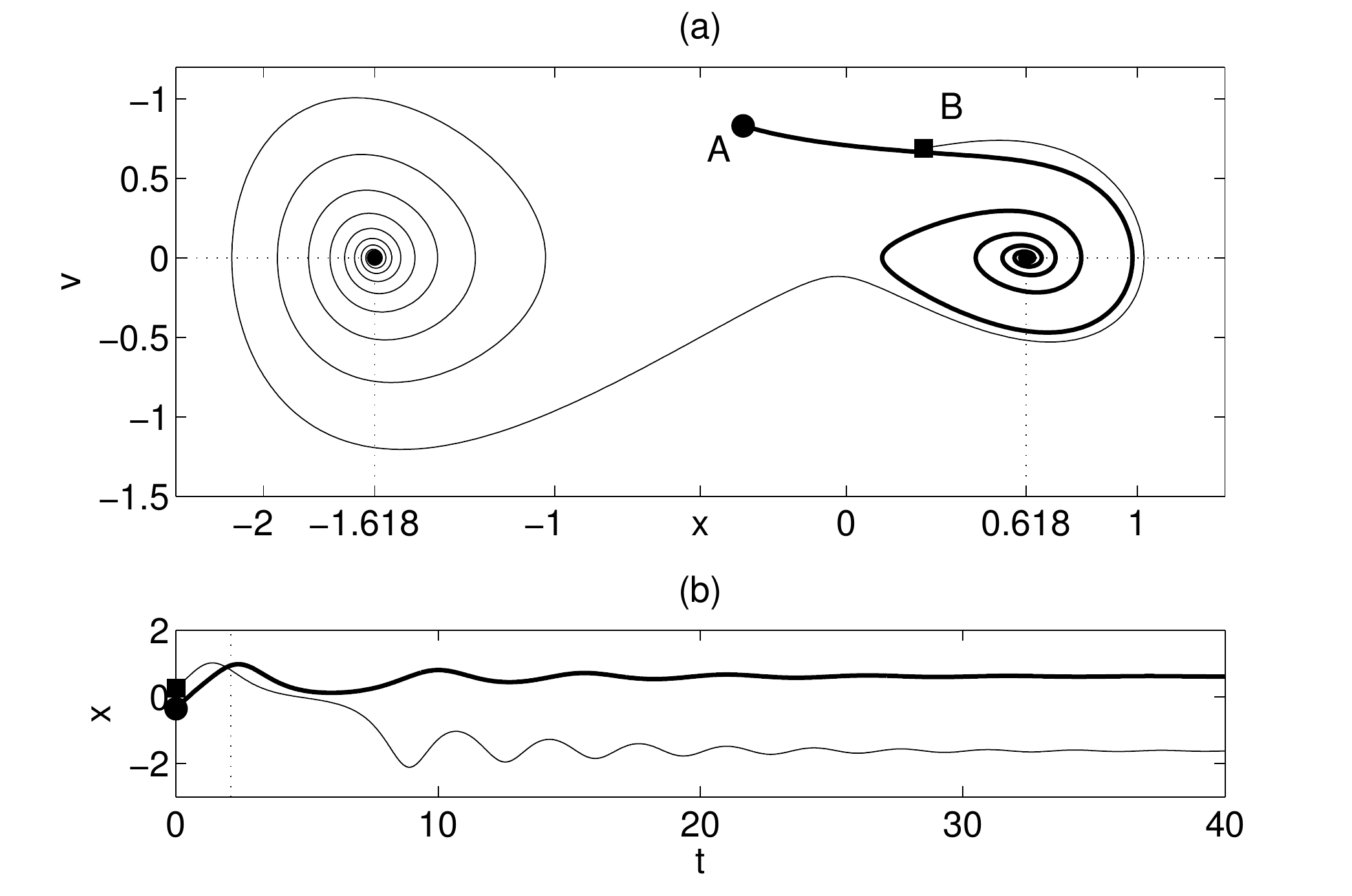}
\caption{Switching-off case: (a) phase portrait and (b) time histories of the two trajectories having points A and B of fig.~\ref{fig:OFFbas}c as initial conditions.}
\label{fig:OFFtime}       
\end{figure}
The same fig.~\ref{fig:OFFbas}b shows in grey the basin of attraction of the forced system (which is the safe zone of the initial system): 
there are points of the basin of attraction of the forced system that are not safe in the transient case and, vice versa, points out of the basin of the forced system which would be safe with the transient switching-off force, as shown in fig.~\ref{fig:OFFbas}b-c and fig.~\ref{fig:OFFtime}.
%
Two relevant points are chosen in order to show the effect of the transient force during the switching-off procedure. 
Point $A$ (black circle) is inside the region of plane contained within the transformed boundary but outside the basin of attraction of the forced system. 
Conversely, point $B$ (black square) lies outside the region limited by the transformed boundary but inside the basin of attraction of the forced system.
On the one hand, if the periodic forcing were not switched off (with no transient procedure), point $A$ would be attracted by the unacceptable periodic attractor because it lies outside the basin of attraction of the forced system. On the other hand, point $B$ would be attracted by the acceptable periodic attractor. 
This is not the case if the effect of the transient force with time-varying amplitude is considered. 
Results from forward time integration are shown both in the phase plane (fig.~\ref{fig:OFFtime}a) and as time histories (fig.~\ref{fig:OFFtime}b). 
The transient force representing the switching-off procedure, acts up to the vertical dotted lines in fig.~\ref{fig:OFFtime}b. Time integration shows that point $A$ converges to the acceptable fixed point, whereas  point $B$ to the unacceptable fixed point of the final system. 

We emphasise that the switching-on problem is akin to the switching off problem because one can be obtained from the other by simple exchange of the initial and final system. 
%
\subsection{Full transient}
The full-transient force  (fig.~\ref{fig:elcentro}) is taken to be proportional to part of the accelelogram of El Centro earthquake, May 1940 \cite{ELCENTRO}. Both initial and final systems are unforced and we want to know the set of conditions in the initial phase plane (at $t=0$) which under the action of the full-transient force give origin to trajectories ending up into the safe region of the final system. As in the switching-off case, the safe zone around the fixed point is chosen coincident with the basin of attraction, capped with the line $v=1$. The boundary of the safe zone of the unforced system is shown as a black line in fig.~\ref{fig:OFFbas}a. The evolution backward in time of the boundary of the safe region is obtained by 
integrating the following system:
\begin{eqnarray}
\label{eq:dotv}
&\dot{x}=   - v, \\
&\dot{v}=   - \left(a(t_{end}-t) - c v + x - x^2 - x^3)\right),\\ 
&a(t)\;\textrm{is the accelelogram in fig.~\ref{fig:elcentro}},\\
&0 \leq t \leq t_{end},\;t_{end}=3.
\end{eqnarray}
The integration of the above equations provides the transformed boundary shown as a dotted line in fig.~\ref{fig:FullBas}a-b.
Fig.~\ref{fig:FullBas}a-b shows in grey the basin of attraction of the unforced system overlapping the transformed boundary.
%
Also in the present case two relevant points are chosen in order to show the effect of the full-transient force. 
Results from forward time integration are shown both in the phase plane (fig.~\ref{fig:FullTime}a) and as time histories (fig.~\ref{fig:FullTime}b). 
The transient force representing an accelelogram, acts up to the vertical dotted lines in fig.~\ref{fig:FullTime}b. The time integration shows that point $A$ converges to the acceptable fixed point, whereas point $B$ to the unacceptable fixed point. 
\begin{figure}
\centering
  \includegraphics[scale=0.6]{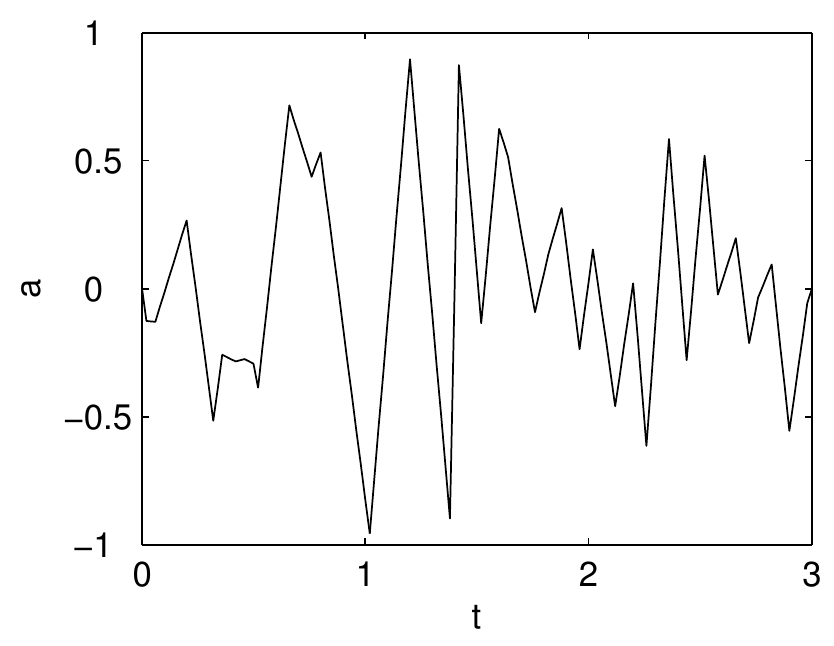}
\caption{Part of the accelelogram of El Centro earthquake, taken from \cite{ELCENTRO}.}
\label{fig:elcentro}       
\end{figure}
\begin{figure}
  \includegraphics[scale=0.425]{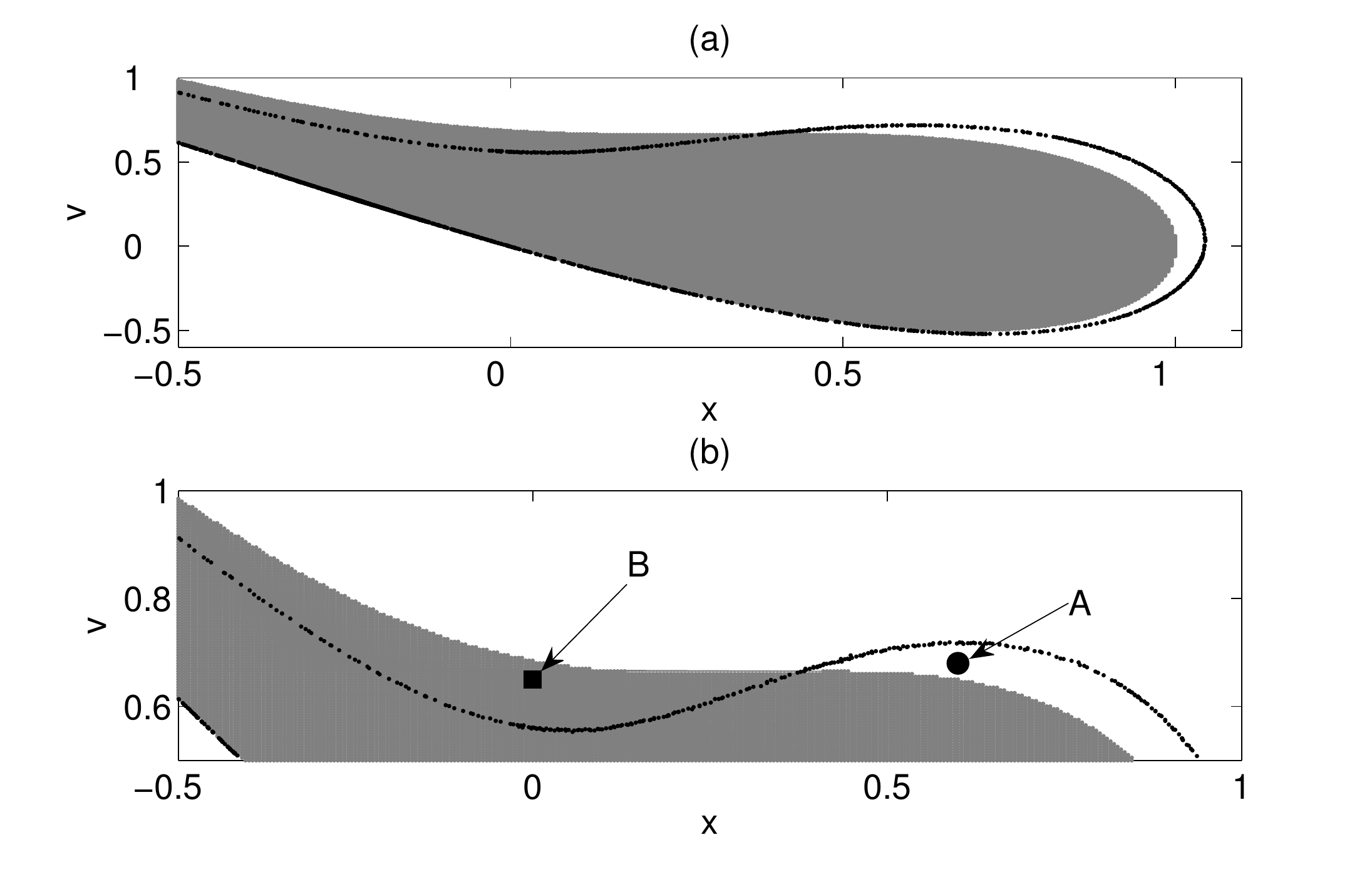}
\caption{Full-transient case: (a) superposition of the transformed boundary of the safe zone of the final system (dotted line) with the basin of attraction of the initial system (grey zone); (b) magnified view of a portion of frame (a).}
\label{fig:FullBas}       
\end{figure}
\begin{figure}
  \includegraphics[scale=0.425]{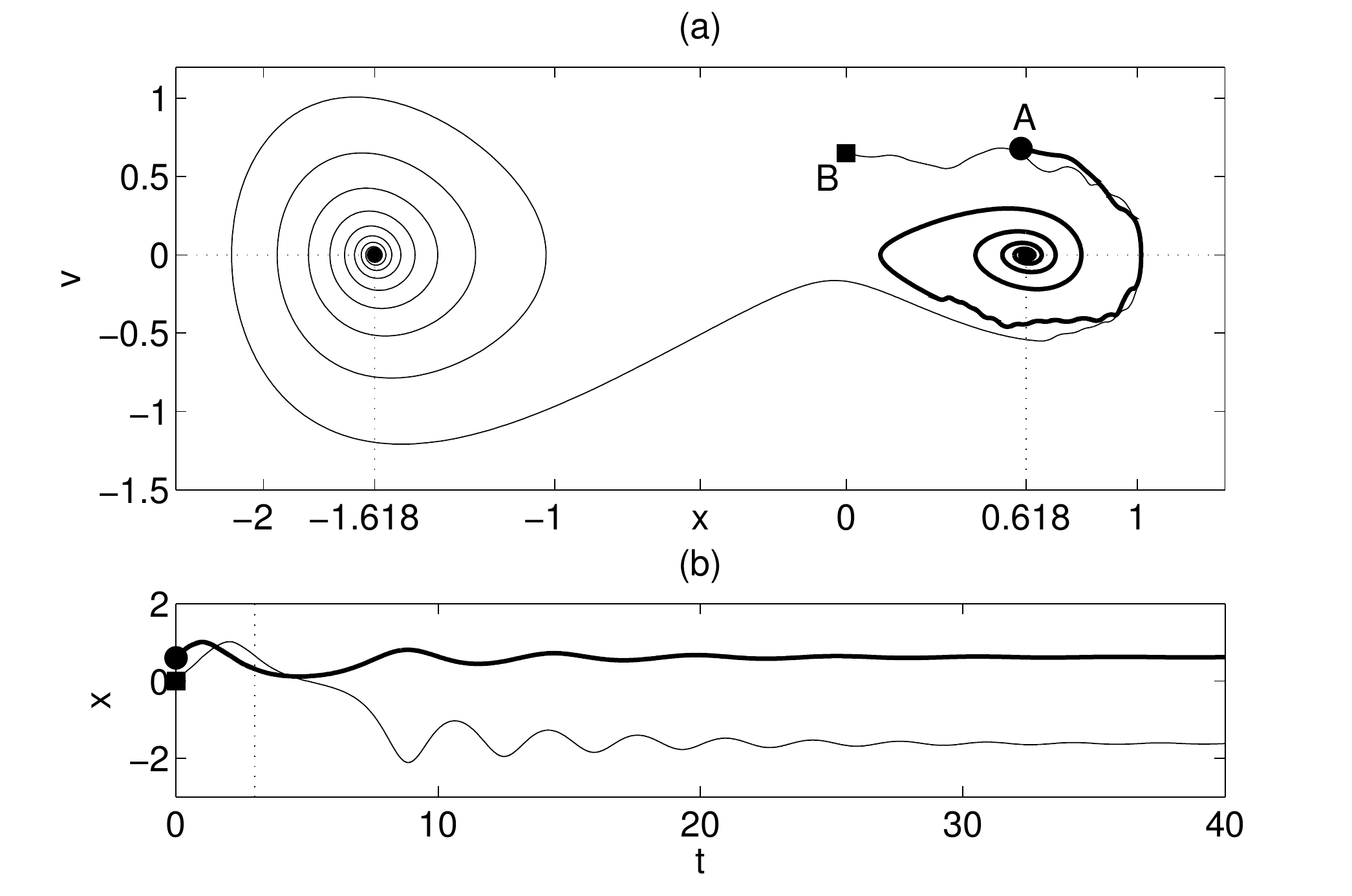}
\caption{Full-transient case: (a) phase portrait and (b) time histories of the two trajectories having points A and B of figure \ref{fig:FullBas}b  as initial conditions.}
\label{fig:FullTime}       
\end{figure}
%
%
%
\section{Conclusions} \label{concl}

This paper presents a preliminary idea to tackle transient problems by using some tools and concepts of the theory of dynamical systems, such as the stable manifolds of saddle points/cycles. 
In particular, we consider systems that are transient because subjected to transient forces. 
A transient force drives a system from a known initial steady state to a target final steady state in a finite time. 
The types of forces considered in this work are representative of many engineering applications, which are switching-off forces and full-transient forces. 
In the switching-off problem, the transient force drives a forced system to an unforced condition. 
In the full-transient problem, the initial and final systems are the same but a transient force acts on the system. 

The key-point is to interpret the action of a transient force as that of a map acting between the two steady systems.
The safe zone, e.g. part of the basin of attraction of a desired attractor, of the target system is mapped backward in time by integration. Doing so, we determine the set of initial conditions which, under the action of the transient force, give origin to a trajectory ending up in the safe zone  of the target final steady system. 
This method has been tested out on a simple mathematical model: a one-degree-of-freedom nonlinear oscillator.

The idea presented in the paper has been explained with the use of diagrams such as those of figures \ref{fig:3}, \ref{fig:OFFbas}  and \ref{fig:FullBas}    
which assume the system to be three-dimensional with states: $x$, $v$ and $t$. 
If the dimension of the system is larger we lose the possibility to visualise the entire safe zones. The detection of the intersection would become a trickier problem, as the number of 
states increases. Moreover, the boundary of the safe zones would be of larger dimension as well. Therefore, the integration backward in time should be applied to a larger number of points, making the procedure more time-consuming.
The method, however, can be applied to many transient systems, as defined by the authors, once one knows: 
i) the duration of the transient;
ii) the mathematical law that describes the transient system;
iii) the initial and final steady systems.

This preliminary work leaves open a few questions. 
First, can a method similar to that presented in this paper be applied in the case of non-deterministic forces, as it is common, for example, in earthquake engineering? 
Secondly, the proposed method could be used to evaluate the size of the portion of the safe zone of the initial system that does not overlap with the image at $t=0$ of the safe zone of the final system. Can this information be used to judge how dangerous different transient forces are? Can a classification of transient forces be introduced, accordingly?
Thirdly, in the examples shown in \S~\ref{numerical} the transformed boundary is a smooth curve and pairs of points which are close in the safe boundary at $t_{end}$ are mapped in pairs of points which are still close at $t=0$. We wonder whether  particular choices of the transient force can exhibit sensitive dependence on initial conditions, or can generate fractal boundaries.
These questions, and possibly many others, will be addressed in future research works.}
%
\begin{acknowledgements}
U.G. would like to acknowledge the financial support from the Italian Ministry of Education, Universities and Research (MIUR) under the PRIN  program 2010/11 N. 2010MBJK5B, Dynamics, Stability and Control of Flexible Structures.

L.M. would like to thank Giulio Ghirardo (Department of Engineering, University of Cambridge) for his help with vector graphics. 

Finally, both authors would like to thank the anonymous reviewers for helpful comments.
\end{acknowledgements}
\end{document}